\documentclass[letterpaper,final]{siamltex}

\usepackage[dvips,dvipdfm]{graphics}
\usepackage{amsfonts}
\usepackage{amssymb}
\usepackage[latin1]{inputenc}
\usepackage{epic, eepic}
\usepackage{url}
\usepackage{pstricks,pst-node,psfrag}
\usepackage{amsfonts,amsmath}
\usepackage{pstricks,pst-node,psfrag}
\usepackage{graphicx}
\usepackage{cite}
\usepackage{url}

\newcommand{\pardiso}   {\textsc{Pardiso{}}}
\newcommand{\metis}     {\textsc{Metis{}}}
\newcommand{\amd}     {\textsc{Amd{}}}

\newcommand{\ilupack}       {\textsc{Ilupack{}}}
\newcommand{\ilupackmatch}       {\textsc{Ilupack-Symmatch{}}}
\newcommand{\cwi}       {\textsc{Cwi{}}}
\newcommand{\jdbsym}       {\textsc{Jdbsym{}}}
\newcommand{\arpack}       {\textsc{Arpack{}}}
\newcommand{\jd}       {\textsc{Jacobi-Davidson{}}}
\newcommand{\nl}{\scriptscriptstyle{0}}
\newcommand{\onexone}   {$1 \times 1$}
\newcommand{\twoxtwo}   {$2 \times 2$}
\newcommand{\PC}{P_{\mathcal{C}}}
\newcommand{\PM}{P_{\mathcal{M}}}
\newcommand{\PMS}{P_{\mathcal{S}}}
\newcommand{\PF}{P_{\rm Fill}}

\newcommand{\includeeps}[2]{\resizebox{#1}{!}{\includegraphics{#2}}}
\newcommand{\figwidths}{0.95\textwidth}
\newcommand{\figwidthd}{0.45\textwidth}

\newcounter{algcnt}
\renewcommand{\thealgcnt}{\thesection.\arabic{algcnt}}

{\refstepcounter{algcnt}{\scshape Algorithm} \thealgcnt. \sl #1}%
{}

\title{On Large Scale Diagonalization Techniques For The
       Anderson Model Of Localization\thanks{This work was supported
       by the Swiss Commission for Technology and Innovation (CTI) under
       grant number 7036 ENS-ES.}}

\author{
    Olaf Schenk\thanks{Department of Computer Science Department,
    University Basel, Klingelbergstrasse 50, CH-4056 Basel, Switzerland,
    ({\tt olaf.schenk@unibas.ch}).}
    \and
    Matthias Bollh\"ofer\thanks{Department of Mathematics, MA 4-5,  TU
    Berlin, Str.\ des 17.\ Juni, 10623 Berlin, Germany,
    ({\tt bolle@math.tu-berlin.de}). Supported by the DFG research center \textsc{Matheon} ``Mathematics for Key Technologies'' in Berlin.}
    \and
    Rudolf A.\ R\"omer\thanks{Centre for Scientific Computing and Department
    of Physics, University of Warwick, Coventry CV4 9BU UK,
    ({\tt r.roemer@warwick.ac.uk}).}
}

\begin{document}


\maketitle \setlength{\doublerulesep}{\arrayrulewidth}

\begin{abstract}
  We propose efficient preconditioning algorithms for an eigenvalue
  problem arising in quantum physics, namely the computation of a few
  interior eigenvalues and their associated eigenvectors for the largest
  sparse real and symmetric indefinite matrices of the Anderson model of
  localization. We compare the Lanczos algorithm in the 1987
  implementation by Cullum and Willoughby with the shift-and-invert
  techniques in the implicitly restarted Lanczos method and in the
  Jacobi-Davidson method. Our preconditioning approaches for the
  shift-and-invert symmetric indefinite linear system are based on
  maximum weighted matchings and algebraic multilevel incomplete $LDL^T$
  factorizations.  These techniques can be seen as a complement to the
  alternative idea of using more complete pivoting techniques for the
  highly ill-conditioned symmetric indefinite Anderson matrices. We
  demonstrate the effectiveness and the numerical accuracy of these
  algorithms. Our numerical examples reveal that recent algebraic
  multilevel preconditioning solvers can accelerative the computation of
  a large-scale eigenvalue problem corresponding to the Anderson model
  of localization by several orders of magnitude.
\end{abstract}



\pagestyle{myheadings}
\thispagestyle{plain}
\markboth{O.~SCHENK, M.~BOLLH\"OFER, AND R.A.R\"OMER}
{LARGE-SCALE DIAGONALIZATION TECHNIQUES \quad $Revision: 1.73 $ \today} 


\section{Introduction}
\label{sect:intro}

One of the hardest challenges in modern eigenvalue computation is the
numerical solution of large-scale eigenvalue problems, in particular
those arising from quantum physics such as, e.g., the Anderson model of
localization (see Section \ref{sect:anderson} for details).  Typically,
these problems require the computation of some eigenvalues and -vectors
for systems which have up to several million unknowns due to their high
spatial dimensions. Furthermore, their underlying structure involves
random perturbations of matrix elements which invalidates simple
preconditioning appraoches based on the graph of the matrices. Moreover,
one is often interested in finding some eigenvalues and associated
eigenvectors in the interior of the spectrum.
The classical Lanczos approach \cite{Par80} has lead to eigenvalue
algorithms \cite{CulW85a,CulW85b} that are in principle able to compute
these eigenvalues using only a small amount of memory. More recent work
on implicitly started Lanczos techniques \cite{LehSY98,arpack} has
accelerated these methods significantly, yet to be fast one needs to
combine this approach with shift-and-invert techniques, i.~e.\ in every
step one has to solve a shifted system of type $A-\sigma I$, where
$\sigma$ is a shift near the desired eigenvalues and
$A\in\mathbb{R}^{n,n}$, $A=A^T$ is the associated matrix. In general
shift-and-invert techniques converge rather quickly which is inline with
the theory \cite{Par80}. Still, an efficient solver is required to solve
systems $(A-\sigma I)x=b$ efficiently with respect to time and memory.
While implicitly restarted Lanczos techniques \cite{LehSY98,arpack}
usually require to solve the system $(A-\sigma I)x=b$ to maximum
precision and thus are mainly suited for sparse direct solvers, the \jd\
method has become an attractive alternative \cite{SleV96} in particular
when dealing with preconditioning methods for linear systems.

Until recently, sparse symmetric indefinite direct solvers were still
far off from symmetric positive definite solvers and this might have
been one major reason why shift-and-invert techniques were not able to
compete with traditional Lanczos techniques \cite{ElsMMR99}, in
particular because of memory constraints. With the invention of fast
matchings-based algorithms \cite{olschowka:1996} which improve the
diagonal dominance of linear systems the situation has dramatically
changed and the impact on preconditioning methods \cite{benzi:2000:phi}
as well as the benefits for sparse direct solvers \cite{sg:04-fgcs} has
been recognized. Furthermore, these techniques have been successfully
transferred to the symmetric case \cite{DufG02,DufP04} allowing modern
state-of-the-art direct solvers \cite{scga:04a} to be orders of
magnitudes faster and more memory efficient than ever, finally leading
to symmetric indefinite sparse direct solvers that are almost as
efficient as their symmetric positive definite counter parts. Recently
this approach has also been utilized to construct incomplete
factorizations \cite{HagS04} with similarly dramatic success.
For a detailed survey on preconditioning techniques for large
symmetric indefinite linear systems the interested reader should
consult \cite{benzi:survey,BenGL05}.




\section{Numerical approach for large systems} \label{sect:contrib}

In the present paper we combine the above mentioned advances with
inverse-based preconditioning techniques \cite{BolS05}. This allows us
to find interior eigenvalues and -vectors for the Anderson problem
several orders of magnitudes faster than traditional algorithms
\cite{CulW85a,CulW85b} while still keeping the amount of memory
reasonably small.

Let us briefly outline our strategy. We will consider recent novel
approaches in preconditioning methods for symmetric indefinite linear
systems and eigenvalue problems and apply them to the Anderson model.
Since the Anderson model is a large-scale sparse eigenvalue problem in
three spatial dimensions, the eigenvalue solvers we deal with are
designed to compute only a few interior eigenvalues and eigenvectors
avoiding a complete factorization. In particular we will use two modern
eigenvalue solvers which we will briefly introduce in Section
\ref{sect:modern}. The first one is \arpack \cite{arpack}, which is a
Lanczos-type method using implicit restarts (cf.\ section
\ref{sect:arnoldi}). We use this algorithm together with a
shift-and-invert technique, i.~e.\ eigenvalues and eigenvectors of
$(A-\sigma I)^{-1}$ are computed instead of those of $A$. \arpack\ is
used in conjunction with a direct factorization method and a multilevel
incomplete factorization method for the shift-and-invert technique.

Firstly, we use the shift-and-invert technique with the novel
symmetric indefinite sparse direct solver that is
part of \pardiso \cite{scga:04a} and we report extensive
numerical results on the  performance of this method.
Section \ref{sect:linear} will give a short overview of the main
concepts that form the \pardiso\ solver.
Secondly, we use \arpack\ in combination with the
multilevel incomplete $LU$
factorization package \ilupack \cite{ilupack2.0}. Here we
present a new {\em indefinite} version of this preconditioner that is
devoted to symmetrically indefinite problems and combines two
basic ideas, namely (i) symmetric maximum weighted matchings
\cite{DufG02,DufP04} and (ii) inverse-based decomposition
techniques \cite{BolS05}. These will be described in Sections \ref{sect:match}
and \ref{sect:iterative}.


As a second eigenvalue solver we use the symmetric version of the
\jd\ method, in particular the implementation \jdbsym
\cite{jdbsym}. This Newton-type method (see section \ref{sect:jd})
is used together with \ilupack \cite{ilupack2.0}.
As we will see in several
further numerical experiments, the synergy of both approaches will
form an extremely efficient preconditioner for the Anderson model
that is memory efficient while at the same time accelerates the
eigenvalue computations significantly: system sizes that resulted
in weeks of computing time \cite{ElsMMR99} can now be computed
within an hour.

\section{The Anderson model of localization}
\label{sect:anderson}

The Anderson model of localization is a paradigmatic model
describing the electronic transport properties of disordered
quantum systems \cite{RomS03,KraM93}. It has been used
successfully in amorphous materials such as alloys \cite{PlyRS03},
semiconductors and even DNA \cite{PorCD04}.  Its hallmark is the
prediction of a spatial confinement of the electronic motion upon
increasing the disorder --- the so-called Anderson {\em
localization} \cite{And58}. When the model is used in $3$ spatial
dimensions, it exhibits a metal-insulator transition in which the
disorder strength $w$ mediates a change of transport properties
from metallic behavior at small $w$ via critical behavior at the
transition $w_{\rm c}$ and on to insulating behavior and strong
localization at larger $w$ \cite{RomS03,KraM93}.
Mathematically, the quantum problem corresponds to a Hamilton
operator in the form of a real symmetric matrix $A$, with quantum
mechanical energy levels given by the eigenvalues $\{\lambda\}$,
and the respective wave functions are simply the eigenvectors of
$A$, i.e.\ vectors $x$ with real entries. With $N = M\times
M\times M$ sites, the quantum mechanical (stationary)
Schr\"{o}dinger equation is equivalent to the eigenvalue equation
$Ax = \lambda x$, which in site representation reads as
\begin{equation}
  \label{eq:eigenequation}
x_{i-1;j;k} + x_{i+1;j;k} + x_{i;j-1;k} + x_{i;j+1;k} + x_{i;j;k-1} + x_{i;j;k+1} +
\varepsilon_{i;j;k} x_{i;j;k} = \lambda x_{i;j;k}
\end{equation}
with $i, j, k = 1, 2, \ldots, M$ denoting the Cartesian coordinates of a
site. The disorder enters the matrix on the diagonal where the entries
$\varepsilon_{i;j;k}$ correspond to a spatially varying disorder
potential and are selected randomly according to a suitable distribution
\cite{KraBMS90}. Here, we shall use the standard box distribution
$\varepsilon_{i;j;k} \in [-w/2, w/2]$ such that the $w$ parameterizes the
aforementioned disorder strength.  Clearly, the eigenvalues of $A$ then
lie within the interval $[-6-w/2, 6+w/2]$. In most studies of the
induced metal-insulator transition, $w$ ranges from $1$ to $30$
\cite{RomS03}. But these values also depend on whether generalizations
to random off-diagonal elements \cite{SouE81,EilRS98a} --- the so-called
random-hopping problem ---, anisotropies \cite{LiKES97,MilRS00} or other
lattice graphs \cite{SchO91,GriRS98} are being considered.

\begin{figure}
\begin{center}
\includegraphics[width=\figwidthd]{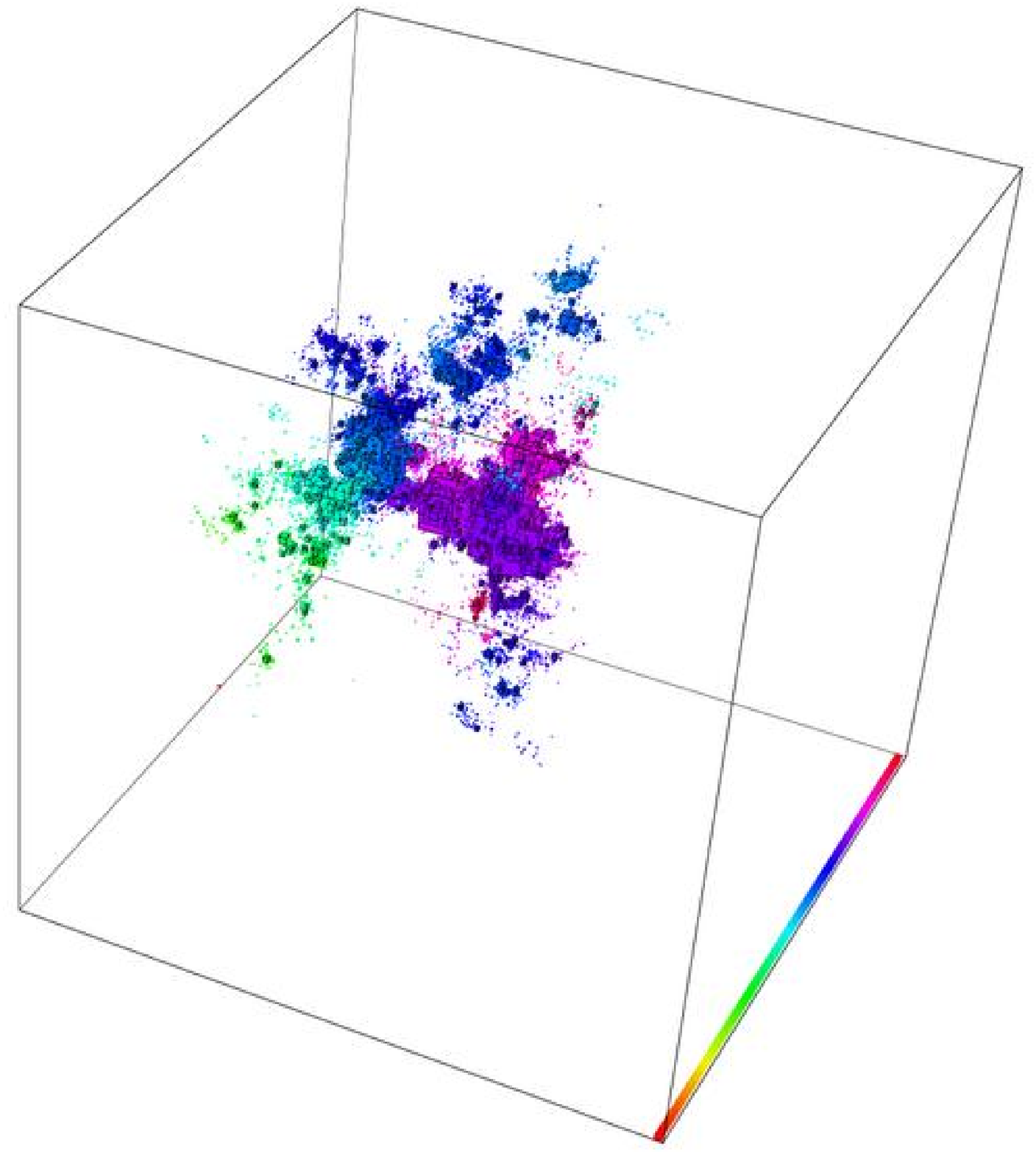}
\includegraphics[width=\figwidthd]{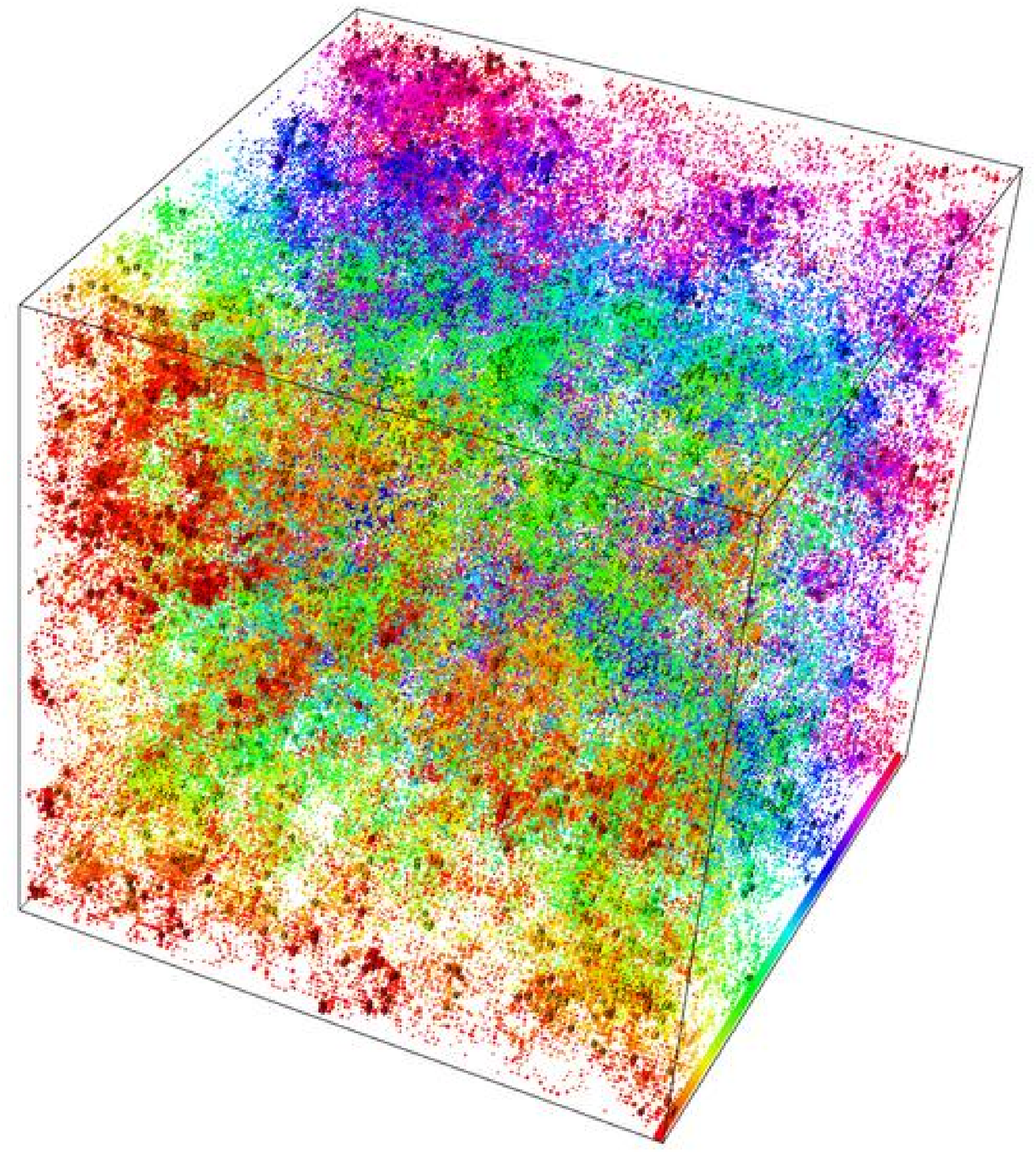}
\end{center}
\caption{\label{fig-extloc_states}
  Extended (left) and localized (right) wave function probabilities for
  the 3D Anderson model with periodic boundary conditions at $\lambda=$
  with $N=100^3$ and $w = 12.0$ and $21.0$, respectively. Every site
  with probability $|x_j|^2$ larger than the average $1/N^3$ is shown as
  a box with volume $|x_j|^2 N$.  Boxes with $|x_j|^2 N>\sqrt{1000}$ are
  plotted with black edges. The color scale distinguishes between
  different slices of the system along the axis into the page.}
\end{figure}

\begin{figure}
\begin{center}
\includegraphics[width=\figwidths]{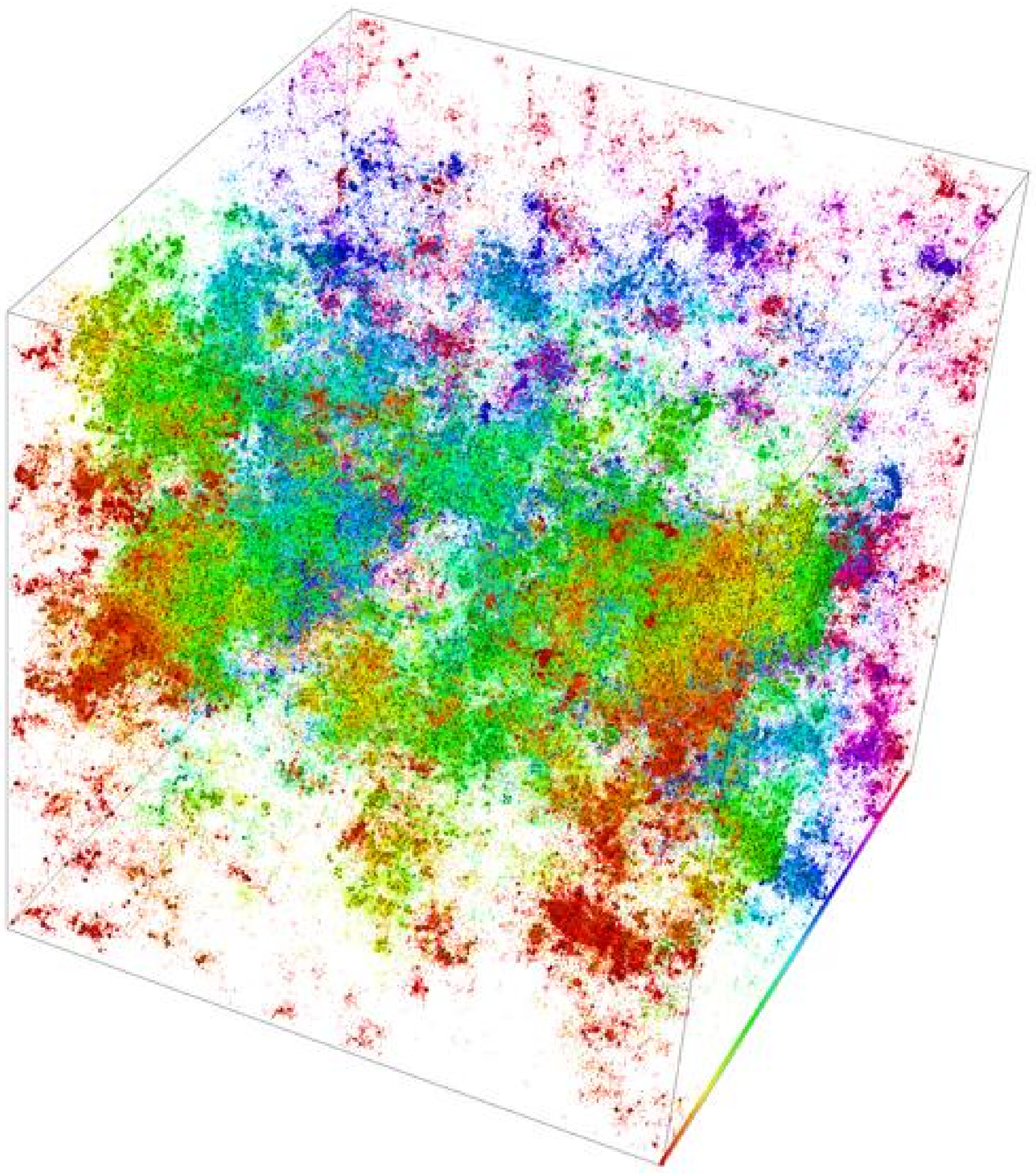}
\end{center}
\caption{\label{fig-crit_state}
  Plot of the electronic eigenstate at the metal-insulator transition
  with $E=0$, $w=16.5$ and $N=250^3$. The box-and-color scheme is as in
  Fig.\ \ref{fig-extloc_states}. Note how the state {\em extends} nearly
  everywhere while at the same time exhibiting certain {\em localized}
  regions of higher $|x_j|^2$ values.}
\end{figure}

The intrinsic physics of the model is quite rich.  For disorders
$w\ll 16.5$, the eigenvectors are extended, i.e.\ $x_{i;j;k}$ is
fluctuating from site to site, but the envelope $|x|$ is
approximately a nonzero constant. For large disorders $w> 16.5$,
all eigenvectors are localized such that the envelope $|x_n|$ of
the $n$th eigenstate may be approximately written as
$\exp-[\vec{r}-\vec{r}_{n}]/l_n(w)$ with $\vec{r} = (i, j, k)^T$
and $l_n(w)$ denoting the {\em localization length} of the
eigenstate. In Figure \ref{fig-extloc_states}, we show examples of
such states. Note that $|x|^2$ and not $x$ corresponds to a
physically measurable quantity and is therefore the observable of
interest to physicists.
Directly at $w = w_{\rm c}\approx 16.5$, the extended states at
$\lambda= 0$ vanish and no current can flow. The wave function vector
$x$ appears simultaneously extended and localized as shown in Fig.\
\ref{fig-crit_state}.

In order to numerically distinguish these three regimes, namely,
localized, critical, and extended behaviors, one needs to (i) go to
extremely large system sizes of order $10^{6}$ to $10^{8}$
and (ii) average over many different
realizations of the disorder, i.e., compute eigenvalues or eigenvectors
for many matrices with different diagonals. In the present paper, we
concentrate on the computation of a few eigenvalues and corresponding
eigenvectors for the physically most interesting case of critical
disorder $w_{\rm c}$ and in the center of $\sigma(A)$, i.e., at
$\lambda= 0$, for large system sizes
\cite{SouE84,Aok86,BraHS96,MilRS97}. Since there is a high density of
states for $\sigma(A)$ at $\lambda = 0$ in all cases, we have the
further numerical challenge of clearly distinguishing the eigenstates in
this high density region.
%
%

\subsection{Lanczos algorithm and the Cullum-Willoughby implementation}
\label{sect:cwi}

Since the mid-eighties, the preferred numerical tool to study the
Anderson matrix and to compute a selected set of eigenvectors,
e.g.\ as needed for a multifractal analysis at the transition, was
provided by the Cullum-Willoughby implementation (\cwi)
\cite{CulW85a,CulW85b,DayTV04} of the Lanczos algorithm.
The algorithm iteratively generates a sequence of orthogonal vectors
$v_i$, $i = 1, \ldots K$, such that $V^T_K A V_K = T_K$, with $V =
[v_1, v_2, \ldots v_K]$ and $T_K$ a symmetric tridiagonal $K\times K$
matrix. In exact arithmetic, the recursion
\begin{equation}
  \label{eq:lanczos}
  \beta_{i+1}v_{i+1} = A v_i - \alpha_i v_i - \beta_i v_{i-1},
\end{equation}
where $\alpha_i = v^T_i A v_i$ and $\beta_{i+1} = v_{i+1} A v_i$ are the
diagonal and subdiagonal entries of $T_K$ and $v_0 = 0$ and $v_1$ is an
arbitrary starting vector, is an orthogonal transformation to
tridiagonal form that needs $K=N$ matrix-vector multiplications. The
eigenvalues of the tridiagonal matrix $T_K$ (Ritz values) are then
simply the eigenvalues of the matrix $A$ and the associated Ritz vectors
yield the eigenvectors.

Since $A$ is sparse and symmetric, the underlying matrix-vector
multiplication on the \cwi\ can be programmed very efficiently,
either by directly coding or appropriate sparse storage schemes
--- only the diagonal needs to be stored in any case.
Additionally, the \cwi\ is a Lanczos implementation in which {\em
no reorthogonalization} is performed. Rather, spurious eigenvalues
are identified by extending the set of Ritz vectors to more than
the $N$ present in exact arithmetic. Typically, we find that $K
\approx 4N$ is sufficient for all ``good'' eigenvalues to have
replicated themselves at least twice --- a further sign that the
algorithm has converged. Clearly, this strategy work well in the
present case since the disorder destroys any symmetry-induced
degeneracies.

The \cwi\ is memory efficient and does not need elaborate
reorthogonalization schemes, but does need to construct many Ritz
vectors which is computationally intensive. Nevertheless, in 1999
\cwi\ was still significantly
faster than more modern iterative schemes
\cite{ElsMMR99}. The main reason for this surprising result lies
in the indefiniteness of $A$, which led to severe difficulties
with solvers more accustomed to standard Laplacian-type problems.

\section{Modern approaches for solving symmetric indefinite eigenvalue problems}
\label{sect:modern}


When dealing with eigenvalues near a given real shift $\sigma$,
the Lanczos algorithm \cite{Par80} is usually accelerated when
being applied to the shifted inverse $(A-\sigma I)^{-1}$ instead
of $A$ directly. This approach relies on the availability of a
fast solution method for linear systems of type $(A-\sigma I)x=b$.
However, the limited amount of available memory only allows for a
small number of solution steps and sparse direct solvers also need
to be memory-efficient to turn this approach into a practical
method.

The limited number of Lanczos steps has lead to modern implicitly
restarted methods \cite{Sor92,LehSY98} which ensure that the
information about the desired eigenvalues is inherited when being
restarted. With increasing number of preconditioned iterative
methods for linear systems \cite{Saa03}, Lanczos-type algorithms
have become less attractive mainly because in every iteration step
the systems of type  $(A-\sigma I)x=b$ have to be solved to {\em full}
accuracy in order to avoid false eigenvalues. In contrast to this,
\jd-like methods \cite{SleV96} allow using a crude approximation
of the underlying linear system. From the point of view of linear
solvers as part of the eigenvalue computation, modern direct and
iterative methods need to inherit the symmetric structure $A=A^T$
while remaining both time and memory efficient. Symmetric
matching algorithms \cite{DufG02,DufP04,scga:04a} have
significantly improved these methods.

\subsection{The shift-and-invert mode of the restarted Lanczos method}
\label{sect:arnoldi}

The Lanczos method for real symmetric matrices $A$ near a shift
$\sigma$ is based on computing successively orthonormal vectors
$[v_1,\dots,v_k,v_{k+1}]$ and a tridiagonal $(k+1)\times k$ matrix
\begin{equation}
\underline{T_k}
=\left(
\begin{array}{cccc}
\alpha_1& \beta_1    &\\
\beta_1& \alpha_2    & \ddots \\
   & \ddots & \ddots  & \beta_{k-1}\\
   &        & \beta_{k-1} &  \alpha_k  \\ \hline
   &        &         & \beta_k
\end{array}
\right)\equiv
\left(
\begin{array}{c}
{T_k}\\ \hline \rule{0ex}{2.2ex} \beta_k e_k^T
\end{array}
\right),
\end{equation}
where $e_k$ is the $k$th unit vector in $\mathbb{R}^k$, such that
\begin{equation}
(A-\sigma
I)^{-1}[v_1,\dots,v_k]=[v_1,\dots,v_k,v_{k+1}]\underline{T_k}.
\end{equation}
Since only a limited number of Lanczos vectors $v_1,\dots,v_k$ can
be stored and since this Lanczos sequence also consists of
redundant information about undesired small eigenvalues,
implicitly restarted Lanczos methods have been proposed
\cite{Sor92,LehSY98} that use implicitly shifted $QR$
\cite{GolV96} exploiting the small eigenvalues of $T_k$ to remove them
out of this sequence without ever forming a single matrix vector
multiplication with $(A-\sigma I)^{-1}$. The new transformed
Lanczos sequence
\begin{equation}
(A-\sigma I)^{-1}[\tilde v_1,\dots,\tilde v_l]=[\tilde
v_1,\dots,\tilde v_l,\tilde v_{l+1}]\underline{\tilde T_l}
\end{equation}
with $l\ll k$ then allows to compute further $k-l$ approximations.
This approach is at the heart of the symmetric version of \arpack
\cite{arpack,LehSY98}.

\subsection{The symmetric \jd\ method}
\label{sect:jd}

One of the major drawbacks of shift-and-invert Lanczos algorithms
is the fact that the multiplication with $(A-\sigma I)^{-1}$
requires solving a linear system to full accuracy. In
contrast to this, \jd-like algorithms \cite{SleV96} are based on a
Newton-like approach to solve the eigenvalue problem. Like the
Lanczos method the search space is expanded step by step solving
the correction equation
\begin{equation}\label{correction}
(I-uu^T) (A-\theta I) (I-uu^T) z = -r \quad \mbox{ such that }
\quad z=(I-uu^T) z
\end{equation}
where $(u,\theta)$ is the given approximate eigenpair and
$r=Au-\theta u$ is the associated residual. Then the search space based
on $V_k=[v_1,\dots,v_k]$ is expanded by reorthogonalizing $z$ with respect
to $[v_1,\dots,v_k]$ and a new approximate eigenpair is computed
from the Ritz approximation $[V_k, z]^T A[V_k, z]$.
When computing several right eigenvectors, the projection $I-uu^T$
has to be replaced by $I-[Q,u][Q, u]^T$ using the already computed
approximate eigenvectors $Q$. This ensures that the new approximate
eigenpair is orthogonal to those that have already been computed.

The most important part of the \jd\ approach is to construct an
approximate solution for (\ref{correction}) such that
\begin{equation}
(I-uu^T) K (I-uu^T)c=d  \quad \mbox{ with }\quad u^Tz=0
\end{equation}
and $K\approx A-\theta I$ that allows a fast solution of the
system $Kx=b$. Here, there is a strong need for robust
preconditioning methods that preserve symmetry and efficiently
solve sequences of linear systems with $K$. If $K$ is itself
symmetric and indefinite, then the simplified QMR method
\cite{FreN95,FreJ97} using the preconditioner
$\left(I-\frac{uw^T}{w^Tu}\right) K^{-1}$, where $Kw=u$ and the system matrix
$\left(I-uu^T\right)(A-\theta I)$ can be used as iterative method.
Note that here the accuracy of the solution of (\ref{correction})
is uncritical until the approximate eigenpair converges
\cite{FokSV00}. This fact has been exploited in \jdbsym
\cite{jdbsym,arbenzgeuss04}. For an overview on \jd\ methods for symmetric
matrices see \cite{Geu02}.

\section{On recent algorithms for solving  symmetric indefinite systems of
            equations}
\label{sect:linear}

We now report on recent improvements in solving symmetric
indefinite systems of linear equations that have significantly
changed sparse direct as well as preconditioning methods. One key
role that made these approaches successful is played by
the use of symmetric matchings that we review in Section
\ref{sect:match}.

\subsection{Sparse direct factorization methods}
\label{sect:direct}

For a long time dynamic pivoting has been a central point for
nonsymmetric sparse linear solvers to gain stability. Therefore,
improvements in speeding up direct factorization methods were
limited to the uncertainties that have arisen from using pivoting.
Certainly techniques like the column elimination  tree
\cite{GilN93,DemEGLL99} have been a useful tool to predict the
sparsity pattern despite pivoting. However, in the symmetric case
the situation becomes more complicated since only symmetric
reorderings applied to both, columns and rows, are required and no
a-priori choice of pivots is given. This makes it almost
impossible to predict the elimination tree in a sensible manner
and the use of cache-oriented level-$3$ BLAS \cite{DodL85,DonDHH88}
was impossible.

With the introduction of symmetric maximum weighted matchings
\cite{DufG02} as alternative to complete pivoting it is now possible to
treat symmetric indefinite systems almost similar to symmetric positive
definite systems. This allows the prediction of fill using the
elimination tree \cite{GeoN85} and thus to set up the data structures
that are required to predict dense submatrices (also known as
supernodes). This in turn means that one is able to exploit level-3 BLAS
applied to the supernodes. Consequently, the classical Bunch-Kaufman
pivoting approach \cite{BunK77} need to be performed only inside the
supernodes.

This approach has recently been successfully implemented in the
sparse direct solver \pardiso \cite{scga:04a} and as a major
consequence, this novel approach has improved the sparse
indefinite solver to become almost as efficient as its symmetric
positive analogy.
Certainly for the Anderson problem studied here, \pardiso\ is
about 2 orders of magnitude more efficient than previously used
direct solvers \cite{ElsMMR99}.
We also note that the idea of symmetric weighted matchings can be
carried over to incomplete factorization methods with similar
success \cite{HagS04}.

\subsection{Symmetric weighted matchings as an alternative to complete
            pivoting techniques}
\label{sect:match}

Symmetric weighted matchings \cite{DufG02,DufP04}, which will be
explained in detail in Section \ref{sec:symmetric}, can be viewed as a
preprocessing step that rescales the original matrix and at the same
time improves the block diagonal dominance. By this strategy, all
entries are at most one in modulus and in addition the diagonal blocks
are either $1\times 1$ scalars $a_{ii}$ such that $|a_{ii}|=1$ (in
exceptional cases we will have $a_{ii}=0$) or they are $2\times 2$
blocks
\[
\left(
\begin{array}{cc}
a_{ii} & a_{i,i+1}\\
a_{i+1,i} & a_{i+1,i+1}
\end{array}
\right)
\mbox{ such that } |a_{ii}|, |a_{i+1,i+1}|\leqslant 1 \mbox{ and }
|a_{i+1,i}|=|a_{i,i+1}|=1.
\]
Although this strategy does not necessarily ensure that symmetric
pivoting like in Ref.\ \cite{BunK77} is unnecessary, it is nevertheless
likely to waive dynamic pivoting during the factorization process.
It has been shown in \cite{DufP04} that based on symmetric
weighted matchings the performance of the sparse symmetric
indefinite {\em multifrontal} direct solver MA57 is improved
significantly, although a dynamic pivoting strategy by Duff and
Reid \cite{DufR83} was still present. Recent results in
\cite{scga:04a} have shown that the absence of dynamic pivoting
does not harm the method anymore and that therefore symmetric
weighted matchings can be considered as alternative to complete
pivoting.

\section{Symmetric reorderings to improve the results of pivoting on
restricted subsets}
\label{sect:weighted-match}

In this section we will discuss weighted graph matchings as an
additional preprocessing step. The motivation for weighted matching
approaches is to identify large entries in the coefficient matrix $A$
that, if permuted close to the diagonal, permit the factorization
process to identify more acceptable pivots and proceed with fewer pivot
perturbations.  These methods are based on maximum weighted matchings
$\mathcal{M}$ and improve the quality of the factor in a complementary
way to the alternative idea of using more complete pivoting techniques.
The idea to use a permutation $\PM$ associated with a weighted matching
$\mathcal{M}$ as an approximation of the pivoting order for nonsymmetric
linear systems was firstly introduced by Olschowka and Neumaier
\cite{olschowka:1996} and extended by Duff and Koster \cite{duko:99a} to
the sparse case.
Permuting the rows $A \leftarrow P_\mathcal{M} A$ of the sparse
system to ensure a zero-free diagonal or to maximize the product
of the absolute values of the diagonal entries are techniques that
are now often regularly used for nonsymmetric matrices
\cite{benzi:2000:phi,superlu:04,sg:04-fgcs,oschenkrogu:04}.

\subsection{Matching algorithms for nonsymmetric matrices}

Let $A = (a_{ij}) \in \mathbb{R}^{n \times n}$ be a general matrix.  The
nonzero elements of $A$ define a graph with edges $\mathcal{E} =
\{(i,j): a_{ij} \neq 0 \}$ of ordered pairs of row and column indices. A
subset $\mathcal{M} \subset \mathcal{E}$ is called a matching or a
transversal, if every row index $i$ and every column index $j$ appear
at most once in $\mathcal{M}$. A matching $\mathcal{M}$ is called
perfect if its cardinality is $n$. For a nonsingular matrix at least one
perfect matching exists and can be found with well known algorithms.
With a perfect matching $\mathcal{M}$, it is possible to define a
permutation matrix $\PM = (p_{ij})$ with:
\begin{equation}
  p_{ij} =
  \begin{cases}
    1 & (j,i) \in \mathcal{M}, \\
    0 & \text{otherwise}.
  \end{cases}
\end{equation}
As a consequence, the permutation matrix $\PM A$ has nonzero
elements on its diagonal.  This method only takes the nonzero
structure of the matrix into account.  There are other approaches
which maximize the diagonal values in some sense.  One possibility
is to look for a matrix $\PM$, such that the product of the
diagonal values of $\PM A$ is maximal. In other words, a
permutation $\sigma$ has to be found, which maximizes
\begin{equation}
  \label{eq:1}
  \prod_{i=1}^n |a_{\sigma(i)i}|.
\end{equation}
This maximization problem is solved indirectly. It can be reformulated 
by defining a matrix $C = (c_{ij})$ with
\begin{equation}
  \label{eq:2}
  c_{ij} =
  \begin{cases}
    \log a_i - \log |a_{ij}| & a_{ij} \neq 0 \\
    \infty                     & \text{otherwise},
  \end{cases}
\end{equation}
where $a_i = \max_j |a_{ij}|$, i.e.\ the maximum element in row $i$ of
matrix $A$. A permutation $\sigma$, which minimizes
$ \sum_{i=1}^n c_{\sigma(i)i} $
also maximizes the product~(\ref{eq:1}).

The minimization problem is known as linear sum assignment problem or
bipartite weighted matching problem in combinatorial optimization.  The
problem is solved by a sparse variant of the Kuhn-Munkres algorithm. The
complexity is $O(n^3)$ for full $n \times n$ matrices and $O(n \tau \log
n )$ for sparse matrices with $\tau$ entries. For matrices whose
associated graph fulfills special requirements, this bound can be
reduced further to $O\left(n^\alpha (\tau + n \log n)\right)$ with
$\alpha < 1$. All graphs arising from finite-difference or
finite-element discretizations meet these conditions~\cite{gupta:99}. As
before, we finally get a perfect matching $\mathcal{M}$ that in turn
defines a nonsymmetric permutation $\PM$.

\begin{figure}[tb]
    $A =$
    \begin{minipage}[c]{0.2\linewidth}
        \includeeps{\linewidth}{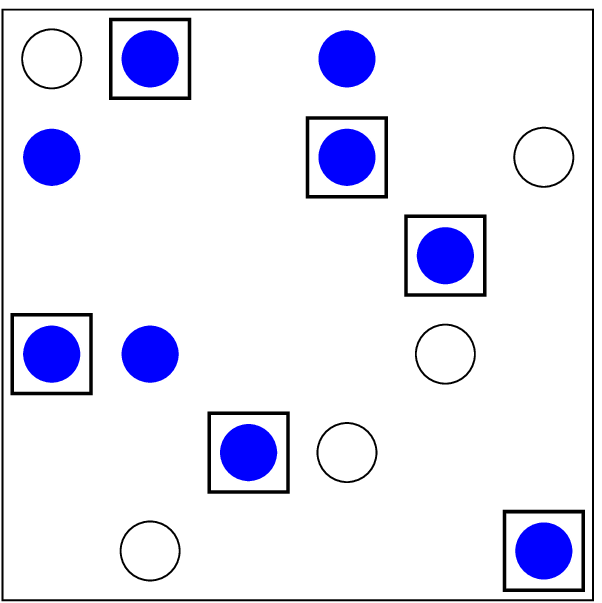}
    \end{minipage}%
    \hfill
    $\PM^T =$
    $\left(
        \begin{array}{cccccc}
        \nl &   1 & \nl & \nl & \nl & \nl \\
        \nl & \nl & \nl &   1 & \nl & \nl \\
        \nl & \nl & \nl & \nl &  1  & \nl \\
          1 & \nl & \nl & \nl & \nl & \nl \\
        \nl & \nl &  1  & \nl & \nl & \nl \\
        \nl & \nl & \nl & \nl & \nl &  1
        \end{array}
    \right) $%
    \hfill
    $\PM A =$
    \begin{minipage}[c]{0.2\linewidth}
        \includeeps{\linewidth}{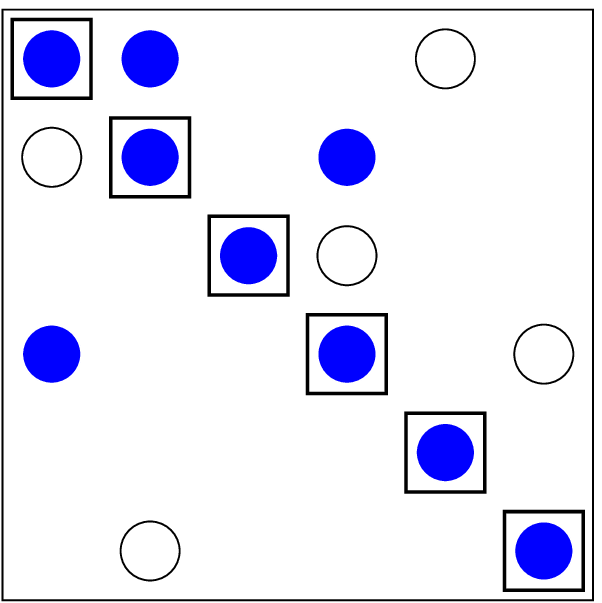}
    \end{minipage}
    \caption{Illustration of the row permutation. A small numerical value
      is indicated by a $\circ$-symbol and a large numerical value by an
      $\bullet$-symbol.  The matched entries $\mathcal{M}$ are marked with squares and
      $\PM=(e_4; e_1; e_5; e_2; e_3; e_6)$.}
    \label{fig:unsym_perm}
\end{figure}

The effect of nonsymmetric row permutations using a permutation associated
with a matching $\mathcal{M}$ is shown in Figure \ref{fig:unsym_perm}. It is
clearly visible that the matrix $\PM A$ is now nonsymmetric, but has the
largest nonzeros on the diagonal.

\subsection{Symmetric \onexone{} and \twoxtwo{} block weighted matchings}
\label{sec:symmetric}

\begin{figure}[tb]
 \begin{center}
    \begin{minipage}[c]{1.0\linewidth}
        \begin{center}
        \begin{minipage}[r]{0.1\linewidth}
         $A$ :
         \end{minipage}
         \begin{minipage}[l]{0.2\linewidth}
           \includeeps{\linewidth}{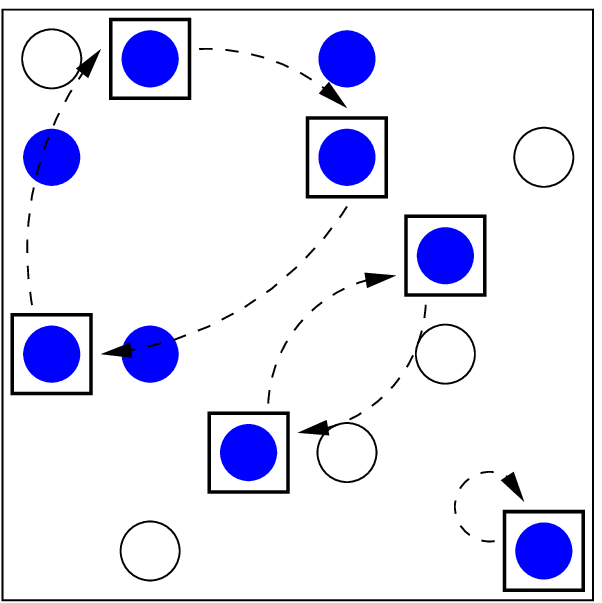}
         \end{minipage}%
    \quad
    \begin{minipage}[r]{0.17\linewidth}
    $\PC A \PC^T$ =
    \end{minipage}
    \begin{minipage}[c]{0.2\linewidth}
        \includeeps{\linewidth}{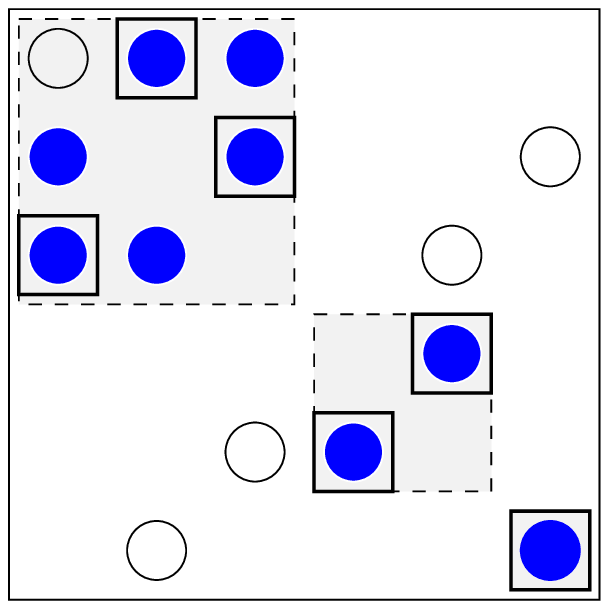}
    \end{minipage}
    \end{center}
    \end{minipage}

    \vspace{0.2cm}

   \begin{minipage}[c]{1.0\linewidth}
       \begin{center}
        \begin{minipage}[r]{0.1\linewidth}
         $A$ :
         \end{minipage}
         \begin{minipage}[l]{0.2\linewidth}
           \includeeps{\linewidth}{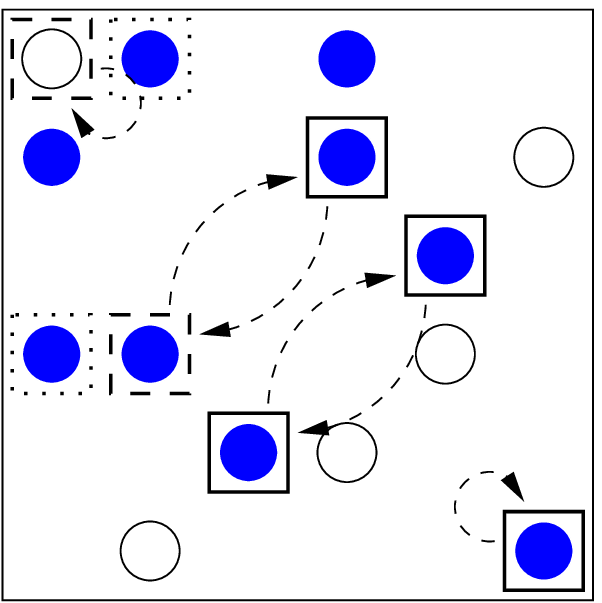}
         \end{minipage}%
    \quad
    \begin{minipage}[r]{0.17\linewidth}
    $\PMS A \PMS^T$ =
    \end{minipage}
    \begin{minipage}[c]{0.2\linewidth}
        \includeeps{\linewidth}{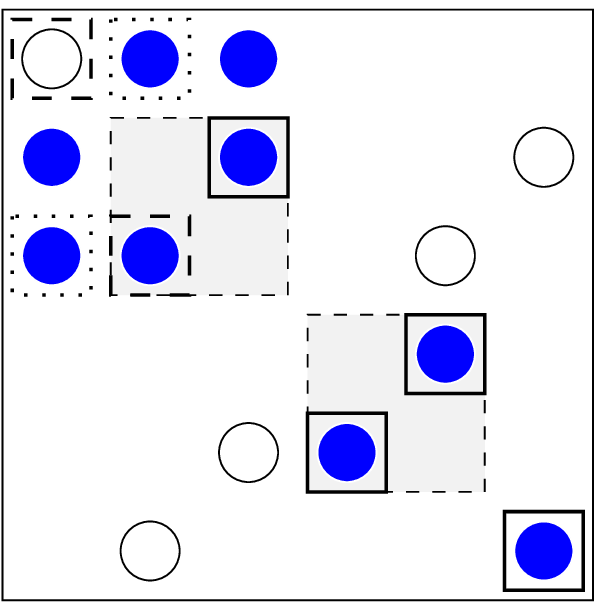}
    \end{minipage}
    \end{center}
    \end{minipage}
 \end{center}

    \caption{Illustration of a cycle permutation with $P_C =
      (e_1;e_2;e_4)(e_3;e_5)(e_6)$ and $\PMS=(e_1)(e_2;e_4)(e_3;e_5)(e_6)$.
      The symmetric matching $ \PMS $
      has two additional elements (indicated by dashed boxes), while one
      element of the original matching fell out (dotted box).  The two
      2-cycles are permuted into \twoxtwo{} diagonal blocks to serve as
      initial  \twoxtwo{} pivots.}
    \label{fig:sym_perm}
\end{figure}

In the case of symmetric indefinite matrices, we are interested in
symmetrically permuting $PAP^T$. The problem is that zero or small diagonal
elements of $A$ remain on the diagonal by using a symmetric permutation $P
AP^T$.  Alternatively, instead of permuting a large\footnote{Large in the
  sense of the weighted matching $\mathcal{M}$.}  off-diagonal element
$a_{ij}$ nonsymmetrically to the diagonal, one can try to devise a
permutation $\PMS$ such that $\PMS A \PMS^T$ permutes this element close to
the diagonal.  As a result, if we form the corresponding \twoxtwo{} block to
$\left[
  \begin{array}{cc}
    a_{ii} &   a_{ij}  \\
    a_{ij} &   a_{jj}  \\
        \end{array}
\right]$, we expect the off-diagonal entry $a_{ij}$ to be large
and thus the \twoxtwo{} block would from a suitable \twoxtwo{}
pivot for the Supernode-Bunch-Kaufman factorization.  An
observation on how to build
    $\PMS$ from the information given by a weighted matching $\mathcal{M}$
    was presented by Duff and Gilbert \cite{DufG02}. They noticed that the
    cycle structure of the permutation $\PM$ associated with the
    nonsymmetric matching $\mathcal{M}$ can be exploited to derive such a
    permutation $\PMS$.  For example, the permutation $\PM$ from Figure
    \ref{fig:unsym_perm} can be written in cycle representation as $ \PC =
    (e_1; e_2; e_4) (e_3; e_5) (e_6)$.  This is shown in the upper graphics
    in Figure \ref{fig:sym_perm}.  The left graphic displays the cycles $(1
    \; 2 \; 4)$, $(3 \;5)$ and $(6)$.  If we modify the original permutation
    $\PM=(e_4; e_1; e_5; e_2; e_3; e_6)$ into this cycle permutation $\PC =
    (e_1; e_2; e_4) (e_3; e_5) (e_6)$ and permute $A$ symmetrically with $
    \PC A \PC^T$, it can be observed that the largest elements are permuted
    to diagonal blocks.  These diagonal blocks are shown by filled boxes in
    the upper right matrix.  Unfortunately, a long cycle would result into a
    large diagonal block and the fill-in of the factor for $\PC A \PC^T$ may
    be prohibitively large. Therefore, long cycles corresponding to $P_{M}$
    must be broken down into disjoint \twoxtwo{} and \onexone{} cycles.
    These smaller cycles are used to define a symmetric permutation $ \PMS =
    (c_1,\dots,c_m)$, where $m$ is the total number of \twoxtwo{} and
    \onexone{} cycles.

    The rule for choosing the \twoxtwo{} and \onexone{} cycles from $\PC$ to
    build $\PMS$ is straightforward.  One has to distinguish between cycles
    of even and odd length.  It is always possible to break down even cycles
    into cycles of length two.  For each even cycle, there are two
    possibilities to break it down. We use a structural metric \cite{DufP04}
    to decide which one to take.  The same metric is also
    used for cycles of odd length, but the situation is slightly different.
    Cycles of length $2l+1$ can be broken down into $l$ cycles of length two
    and one cycle of length one.  There are $2l+1$ different possibilities
    to do this. The resulting $2\times2$ blocks will contain the matched
    elements of $\mathcal{M}$.  However, there is no guarantee that the
    remaining diagonal element corresponding to the cycle of length one will
    be nonzero. Our implementation will randomly select one element as a
    \onexone{} cycle from an odd cycle of length $2l+1$.

A selection of $\PMS$ from a weighted matching $\PM$ is
illustrated in Figure \ref{fig:sym_perm}. The permutation
associated with the weighted matching, which is sorted according
to the cycles, consists of $ \PC = (e_1; e_2; e_4) (e_3; e_5)
(e_6)$. We now split the full cycle of odd length three into two
cycles $(1)(2 4)$ --- resulting in $\PMS = (e_1) (e_2; e_4)
(e_3;e_5) (e_6) $.  If $\PMS$ is symmetrically applied to $A
\leftarrow \PMS A \PMS^T$, we see that the large elements from the
nonsymmetric weighted matching $\mathcal{M}$ will be permuted
close to the diagonal and these elements will have more chances to
form good initial \onexone{} and \twoxtwo{} pivots for the
subsequent (incomplete) factorization.

Good fill-in reducing orderings $\PF$ are equally important for symmetric
indefinite systems. The following section 
introduces two strategies to combine these reorderings with the symmetric
graph matching permutation $\PMS$. This will provide good initial pivots for
the factorization as well as a good fill-in reduction permutation.

\subsection{Combination of orderings $\PF$ for fill reduction with orderings
  $\PMS$ based on weighted matchings}
\label{sect:construct-symperm}

In order to construct the factorization efficiently, care has to
be taken that not too much fill-in is introduced during the
elimination process. We now examine two algorithms for the
combination of a permutation $\PMS$ based on weighted matchings to
improve the numerical quality of the coefficient matrix $A$ with a
fill-in reordering $\PF$ based on a nested dissection from
\textsc{Metis} \cite{karypis:98}.  The first method is based on
compressed subgraphs and has also been used by Duff and Pralet in
\cite{DufP04} in order to find good scalings and orderings for
symmetric indefinite systems.


In order to combine the permutation $\PMS$ with a fill-in reducing
permutation, we compress the graph of the reordered system $\PMS A \PMS^T$
and apply the fill-in reducing reordering to the compressed graph. In the
compression step, the union of the structure of the two rows and columns
corresponding to a \twoxtwo{} diagonal block are built, and used as the
structure of a single, compressed row and column representing the original
ones.

If $G_A=(V;E)$ is the undirected graph of $A$ and a cycle consists
of two vertices $(s,t) \in V$, then graph compression will be done
on the \onexone{} and \twoxtwo{} cycles, which have been found
using a weighted matching $\mathcal{M}$ on the graph.  The
vertices $(s,t)$ are replaced with a single supervertex $\emph{u}
= \{ s, t \} \in V_c$ in the compressed graph $G_c = (V_c,E_c)$.
An edge $e_c = (s, t) \in E_c$ between two supervertices $\emph{s}
= \{ s_1, s_2 \} \in V_c$ and $\emph{t} = \{ t_1, t_2 \}\in V_c$
exists if at least one of the following edges exist in $E: \: (s_1,
t_1), (s_1, t_2), (s_2, t_1)$ or $(s_2, t_2)$. The fill-in
reducing ordering is found by applying \metis {} on the compressed
graph $G_c = (V_c,E_c)$. Expansion of that permutation to the
original numbering yields $P_{fill}$. Hence all \twoxtwo{} cycles
that correspond to a suitable \twoxtwo{} pivot block are reordered
consecutively in the factor.

\section{Symmetric multi-level preconditioning techniques}
\label{sect:iterative}

We now present a new symmetric indefinite approximate multilevel
factorization that is mainly based on three parts which are
repeated in a multilevel framework in each subsystem. The
components consist of (i) reordering of the system, (ii)
approximate factorization using inverse-based pivoting and, (iii)
recursive application to the system of postponed updates.

\subsection{Reordering the given system} The key ingredient
to turn this approach into an efficient multilevel solver consists
of the symmetric maximum weight matching presented in Section
\ref{sect:match}. After the system is reordered into a
representation
\begin{equation}
P_s^T D A D P_s =\hat A,
\end{equation}
where $D, P_s\in\mathbb{R}^{n,n}$, $D$ is a diagonal matrix and
$P_s$ is a permutation matrix, $\hat A$ is expected to have many
diagonal blocks of size $1\times 1$ or $2\times 2$ that are
well-conditioned. Once the diagonal block of size $1\times 1$ and
$2\times 2$ are built, the associated block graph of $\hat A$ is
reordered by a symmetric reordering, e.~g.\ \amd\ \cite{AmeDD96}
or \metis \cite{karypis:98}, i.~e.\
\begin{equation}
\Pi^T P_s^T D A D P_s \Pi =\tilde A,
\end{equation}
where $\Pi\in \mathbb{R}^{n,n}$ refers to the associated symmetric block
permutation.

\subsection{Inverse-based pivoting}

Given $\tilde A$ we compute an incomplete factorization
$LDL^T=\tilde A +E$ of  $\tilde A$. To do this at step $k$ of the
algorithm we have
\begin{equation}
\tilde A=
\left(
\begin{array}{cc}
B & F^T\\
F & C
\end{array}
\right)
=
\left(
\begin{array}{cc}
L_B & 0\\
L_F & I
\end{array}
\right)
\left(
\begin{array}{cc}
D_B & 0\\
  0 & S_C
\end{array}
\right)
\left(
\begin{array}{cc}
L_B^T & L_F^T\\
0& I
\end{array}
\right),
\end{equation}
where $L_{B}\in\mathbb{R}^{k-1,k-1}$ is lower triangular with unit
diagonal and $D_{B}\in\mathbb{R}^{k-1,k-1}$ is block diagonal with
diagonal blocks of size $1\times1$ and $2\times 2$. Also,
$S_C=C-L_FD_BL_F^T=\left(s_{ij}\right)_{i,j}$ denotes the
approximate Schur complement. To proceed with the incomplete
factorization we perform either a $1\times 1$ update or a $2\times
2$ block update. One possible choice could be to use Bunch's
algorithm \cite{Bun74}. This approach has been used in Ref.\
\cite{HagS04}. Here we use a simple criterion based on block
diagonal dominance of the leading block column.
 Depending on the values
\begin{equation}
d_1=\sum_{j>1} \frac{|s_{j1}|}{|s_{11}|}, \quad
d_2=\sum_{j>2}
\|\left(s_{j1}, s_{j2}\right) \left(\begin{array}{cc}
s_{11} & s_{12}\\
s_{12} & s_{22}
\end{array}
\right)^{-1}\|,
\end{equation}
we perform a $2\times 2$ update only if $d_2<d_1$. The leading two
columns of $S_C$ can be efficiently computed using linked lists
\cite{LiSC03} and it is not required to have all entries of $S_C$
available.

When applying the (incomplete) factorization $LDL^T$ to $\tilde A$
we may still encounter a situation where at step $k$ either
$1/|s_{11}|$ or $\|\left(s_{ij}\right)_{i,j\leqslant 2}^{-1}\|$ is
large or even infinity. Since we are dealing with  an incomplete
factorization we propose to use inverse-based pivoting
\cite{BolS05}. Therefore we require in every step that
\begin{equation}
\|\left(
  \begin{array}{cc}
   L_B     & 0\\
   L_E & I
  \end{array}
\right)^{-1}\| \leqslant \kappa
\end{equation}
for a prescribed bound $\kappa$. If after the update using a
$1\times 1$ pivot (or $2\times 2$ pivot) the norm of the inverse
lower triangular factor fails to be less than $\kappa$, the update
is postponed and the leading rows/columns of $L_E$, $S_C$ are
permuted to the end. Otherwise depending on whether a $1\times1$
or a $2\times 2$ pivot has been selected, the entries
\begin{equation}
(s_{j1}/s_{11})_{j>1}, \quad \left(
 \left(s_{j1}, \; s_{j2}\right)
 \left(
  \begin{array}{cc}
   s_{11}&s_{12}\\
   s_{12}& s_{22}
  \end{array}
 \right)^{-1}\right)_{j>2}
\end{equation}
become the next (block) column of $L$ and we drop these entries
whenever their absolute value is less than $\varepsilon/\kappa$
for some threshold $\varepsilon$. For a detailed description see
Ref.\ \cite{BolS05}. The norm of the inverse can cheaply be
estimated using a refined strategy of Ref.\ \cite{CliMSW79} and is
part of the software package \ilupack\ that is now extended to the
symmetric indefinite case \cite{ilupack2.0}.

\subsection{Recursive application}

After the inverse-based ILU we have an approximate factorization
\begin{equation}
Q^T \tilde AQ = \left(
\begin{array}{cc}
L_{11} & 0\\
L_{21} & I
\end{array}
\right)
\left(
\begin{array}{cc}
D_{11} & 0\\
  0 & S_{22}
\end{array}
\right)
\left(
\begin{array}{cc}
L_{11}^T & L_{21}^T\\
0& I
\end{array}
\right)
\end{equation}
and it typically does not pay off to continue the factorization
for the remaining matrix $S_{22}$ which consists of the previously
postponed updates. Thus $S_{22}$ is now explicitly computed and
the strategies  for reordering, scaling and factorization are
recursively applied to $S_{22}$ leading to a multilevel
factorization.

Note that in order to save memory $L_{21}$ is not stored but
implicitly approximated by $\tilde
A_{21}(L_{11}D_{11}U_{11})^{-1}$. In addition we use a technique
called \textsl{aggressive dropping} that sparsifies the triangular
factor $L$ a posteriori. To do this observe that when applying a
perturbed triangular factor $\tilde L^{-1}$ for preconditioning
instead of $L^{-1}$ we have
\[
\tilde{L}^{-1}=(I+E_L)L^{-1}, \mbox{ where } E_L=\tilde{L}^{-1} (L-\tilde{L}).
\]
We can expect that $\tilde{L}^{-1}$ serves as a good approximation
to $L^{-1}$ as long as $\|E_L\|\ll 1$. If we obtain $\tilde{L}$
from $L$ by  dropping some entry, say $l_{ij}$ from $L$, then we
have to ensure that
\[
\|\tilde{L}^{-1}e_i\| \cdot  |l_{ij}|\leqslant \tau \ll 1,
\]
for some moderate constant $\tau<1$, e.g. $\tau=0.1$. To do this
it is required to have a good estimate for
$\nu_i\approx\|\tilde{L}^{-1}e_i\|$ available, for any
$i=1,\dots,n$. In principle it can be computed
\cite{BolS05,CliMSW79} using $\tilde L^\top$ instead of $\tilde
L$. Last, knowing how many entries exist in column $j$, we could
drop any $l_{ij}$ such that
\[
|l_{ij}| \leqslant \tau/(\nu_i \cdot \#\{l_{kj}:\; l_{kj}\not=0, k=j+1,\dots,n\}).
\]

\subsection{Iterative solution}

By construction, the computed incomplete multilevel factorization is symmetric
but indefinite. For the iterative solution of linear systems using the
multilevel factorization, in principle different Krylov subspace solvers
could be used. For example, general
methods that do not explicitly use symmetry
(e.g. GMRES \cite{SaaS86}) or
methods like SYMMLQ \cite{PaiS75} which preserve the symmetry of the orginal
 matrix but which are only devoted for symmetric positive
definite preconditioners.
To fully exploit both, symmetry and indefiniteness at the same time,
here the simplified QMR method \cite{FreN95,FreJ97} is chosen.

\section{Numerical Experiments}
\label{sect:exp}

Here we present numerical experiments that show that the
previously outlined advances in symmetric indefinite sparse direct
solvers as well as in preconditioning methods significantly
accelerate modern eigenvalue solvers  and allow us to gain orders
of magnitude in speed compared to more conventional methods.

\subsection{Computing Environments and Software}

All large scale numerical experiments for the Anderson model of
localization were performed on an SGI Altix 3700/BX2 with 56 Intel
Itanium2 1.6 GHz processors and 112 GB of memory. If not
explicitly stated, we always used only one processor of the system
and all algorithms were implemented in either C or Fortran77.  All
codes were compiled by the Intel V8.1 compiler suite using {\em
ifort} and {\em icc} with the -O3 optimization option and linked
with basic linear algebra subprograms optimized for Intel
architectures.
For completeness, let us recall the main software packages used.
\begin{itemize}
\item {\protect\arpack} is a collection of Fortran77 subroutines
designed to solve large scale eigenvalue problems. The eigenvalue
solver has been developed at the Department of Computational and
Applied Mathematics at Rice University.  It is available at
\url{http://www.caam.rice.edu/software/ARPACK}.

\item \jdbsym\ is a C library implementation of the \jd\ method
optimized for symmetric eigenvalue problems. It solves
eigenproblems of the form $Ax = \lambda x$ and $Ax = \lambda B x$
with or without preconditioning, where A is symmetric and B is
symmetric positive definite. It has been developed at the Computer
Science Department of the ETH Zurich and is available at
\url{http://people.web.psi.ch/geus/software.html}.

\item \pardiso\ is a fast direct solver package, developed at the
Computer Science Department of the University of Basel and
available at
\url{http://www.computational.unibas.ch/cs/scicomp/software/pardiso}.

\item \ilupack\ is an algebraic multilevel preconditioning
software package. This iterative solver has been developed at the
Mathematics Department of the Technical University of Berlin. It
is available at \url{http://www.math.tu-berlin.de/ilupack}.

\end{itemize}

\subsection{Numerical Results}

In our numerical experiments we will first compare the classical
\cwi {} with the shift-and-invert Lanczos method using implicit
restarts. The latter is part of \arpack \cite{arpack}. For the
solution of the symmetric indefinite system $A-\theta I$ we use
the most recent version of sparse direct solver
\pardiso \cite{scga:04a}. This version is based on symmetric
weighted matchings and uses \metis{} as symmetric reordering
strategy. The numerical results deal with the computation of $5$
eigenvalues of the Anderson matrix $A$ near $\lambda=0$. Here we
state the results for the physically most interesting critical
disorder strength $w=w_c=16.5$ (cf.\ Table
\ref{tab:cpu-pardiso-cwi}). As can be seen from Table
\ref{tab:cpu-pardiso-cwi}, the
\pardiso-based shift-and-invert Lanczos is clearly superior to the
classic \cwi {} method by at least one order of magnitude
regarding computation time. Despite this success, with increasing
problem size the amount of memory consumed by the sparse direct
solver becomes significant and numerical results $N$ larger than
$1'000'000$ are skipped.

\begin{table}[htb]
    \centering
    \caption{CPU times in seconds and memory requirements in GB to compute
      at $w=16.5$ five eigenvalues  closest to $\lambda=0$ of an Anderson
      matrix of size $N=M^3 \times M^3$ with \cwi {} and
      \arpack-\pardiso. For \cwi\ and $M=90, 100$, not all $5$ eigenvalues converged successfully, so the eigenvector reconstruction finished quicker, leading to apparently shorter CPU times ($^*$).}
    \begin{tabular}{|r|r|r|r|r|r|}\hline
        \hline
   \multicolumn{1}{|c|}{M}&\multicolumn{1}{|c|}{N}& \multicolumn{2}{c|}{\cwi} &\multicolumn{2}{c|}{\arpack-\pardiso}\\
    &           &   \multicolumn{1}{|c|}{time}   &
   \multicolumn{1}{|c|}{mem.} & \multicolumn{1}{|c|}{time}   & \multicolumn{1}{|c|}{mem.} \\ \hline\hline
 30 &    27'000 &      21   &     0.01 &     9 &   0.08 \\
 40 &    64'000 &     300   &     0.02 &    46 &   0.28 \\
 50 &   125'000 &   1'246   &     0.04 &   157 &   0.68 \\
 60 &   216'000 &   4'748   &     0.07 &   495 &   1.49 \\
 70 &   343'000 &  15'100   &     0.11 & 1'309 &   3.00 \\
 80 &   512'000 &  39'432   &     0.16 & 3'619 &   5.12 \\
 90 &   729'000 & 97'119$^*$&     0.23 & 7'909 &   8.70 \\
100 & 1'000'000 &255'842$^*$&     0.32 &20'239 &  14.34 \\  \hline\hline
 \end{tabular}
    \label{tab:cpu-pardiso-cwi}
\end{table}


Instead, we switch to the \ilupack-based preconditioner that is
also based on symmetric weighted matchings and in addition uses
inverse-based pivoting. In particular, for our experiments we use
$\kappa=5$ as bound for the norm $\|L^{-1}\|$ of the inverse
triangular factor and \amd\ for the symmetric reordering. We also
tried to use \metis{} but for this particular matrix problem we
find that \amd\ is clearly more memory efficient. Next we compare
\pardiso-based shift-and-invert Lanczos (\arpack) with that using
\ilupack\ and the simplified QMR as inner iterative solver. Here
we use $\varepsilon=1/sqrt(N)$ with aggressive dropping and the
QMR method is stopped once the norm of residual satisfies
$\|Ax-b\|\leqslant 10^{-10}\|b\|$. In order to illustrate the
benefits of using symmetric weighted matchings we also tried
\ilupack\ without matching, but the numerical results are
disappointing as can be seen from the $\dagger$s in Table
\ref{tab:all}.  We emphasize that the multi-level approach is
crucial, a simple use of incomplete factorization methods
 without multi-level preconditioning \cite{HagS04} does not give
the desired results. Besides the effect of matchings we also
compare how the performance of the methods changes when varying
the value $w$ from the critical value $w=w_c=16.5$ to $w=12.0$ and
$w=21.0$.  We find that these changes do not affect the sparse
direct solver at all while the multilevel ILU significantly varies
in its performance. Up to now our explanation for this effect is
the observation that with increasing $w$ the diagonal dominance of
the system also increases and the \ilupack\ preconditioner gains
from higher diagonal dominance. As we can see from Table
\ref{tab:all}, \ilupack\ still uses significantly less memory than
the direct solver \pardiso\ for all values of $w$ and it is the
only method we were able to use for larger $N$ due to the memory
constraints. Also, the computation time is best.

\begin{table}[htb]
\centering
\caption{CPU times in seconds  and memory requirements in GB to
         compute five eigenvalues closest to
         $\lambda=0$ of an Anderson matrix of size $N^3 \times N^3$
         with \arpack-\pardiso, \arpack-\ilupack, and
         \arpack-\ilupackmatch. The symbol '---' indicates that a memory
         consumption was larger than $25$ GB and '$\dagger$' indicates
         memory problems with respect to the fill--in.}
\begin{tabular}{|r|c|r|r|r|r|r|r|}\hline
        \hline
   \multicolumn{1}{|c|}{M} & \multicolumn{1}{|c|}{W} &  \multicolumn{6}{c|}{\arpack}\\
      &      &  \multicolumn{2}{c|}{\pardiso} &
                \multicolumn{2}{c|}{\ilupack} &
                \multicolumn{2}{c|}{ \ilupackmatch}\\
      &      &   \multicolumn{1}{|c|}{time}  & \multicolumn{1}{|c|}{mem.} &
                 \multicolumn{1}{|c|}{time}  & \multicolumn{1}{|c|}{mem.} &
                 \multicolumn{1}{|c|}{time}  & \multicolumn{1}{|c|}{mem.}\\ \hline\hline
   70 & 12.0 &    1'359 &   3.00 &       5'117 &     1.09 &      2'140 &    0.95 \\
  100 & 12.0 &   20'639 &  14.34 &      39'222 &     5.62 &     13'583 &    3.20 \\
  130 & 12.0 &    ---   &  ---   &    $\dagger$& $\dagger$&     65'722 &         \\
\hline\hline
   70 & 16.5 &    1'305 &   3.00 &         504 &     0.33 &        477 &   0.31  \\
  100 & 16.5 &   20'439 &  14.34 &       2'349 &     0.95 &      2'177 &   0.89  \\
  130 & 16.5 &   ---    &  ---   &       6'320 &     2.09 &      6'530 &   1.95  \\
  160 & 16.5 &   ---    &  ---   &      23'663 &     3.95 &     13'863 &   3.63  \\
\hline\hline
   70 & 21.0 &    1'225 &   3.00 &         371 &     0.22 &        310 &   0.22  \\
  100 & 21.0 &   20'239 &  14.34 &       1'513 &     0.64 &      1'660 &   0.65  \\
  130 & 21.0 &   ---    &  ---   &       3'725 &     1.41 &      3'527 &   1.44  \\
  160 & 21.0 &   ---    &  ---   &      15'302 &     2.63 &     20'120 &   2.68  \\
\hline\hline
\end{tabular}
\label{tab:all}
\end{table}

When using preconditioning methods inside shift-and-invert Lanczos
we usually have to solve the inner linear system for $A-\theta I$
up to the machine precision to make sure that the eigenvalues and
eigenvectors are sufficiently correct. In contrast to this the
\jd\ method allows to solve the associated correction equation
less accurately and only when convergence takes place a more
accurate solution is required.
In order to show the significant difference between the iterative parts of \arpack\
and \jd\
we state the number of iteration steps in Table \ref{tab:it-arpack-jd}.
If we were to aim for more eigenpairs, we
expect that eventually the \jdbsym\ becomes less efficient and
should again be replaced by \arpack.

\begin{table}[htb]
\centering \caption{Number of inner/outer interaction steps inside
\arpack and \jd. The symbol '---' indicates that the computations
were not performed
         anymore for \arpack.}
\begin{tabular}{|r|c|r|r|r||r|r|r|}\hline
        \hline
   \multicolumn{1}{|c|}{M} & \multicolumn{1}{|c|}{W} &  \multicolumn{6}{c|}{\ilupackmatch}\\
    &  &  \multicolumn{3}{c|}{\arpack}&  \multicolumn{3}{c|}{\jd}\\
      &      &   \multicolumn{1}{|c|}{outer}  & \multicolumn{1}{|c|}{total} &
                 \multicolumn{1}{|c|}{inner}  & \multicolumn{1}{|c|}{outer} &
                 \multicolumn{1}{|c|}{total}  & \multicolumn{1}{|c|}{inner}\\
      &      &    &  &
                 \multicolumn{1}{|c|}{average}  & &
                   & \multicolumn{1}{|c|}{average}\\ \hline\hline
   70 & 12.0 &    42    &    871 &   20.7    &        43 &   246 &    5.7 \\
  100 & 12.0 &    43    &   1101 &   25.6    &        44 &   325 &    7.4 \\
  130 & 12.0 &    42    &   1056 &   25.1    &        44 &   274 &    6.2 \\
\hline\hline
   70 & 16.5 &     43   &   611  &   14.2    &   41      & 200  &  4.88   \\
  100 & 16.5 &     43   &   857  &   19.9    &   43      & 231  &  5.37   \\
  130 & 16.5 &     42   &  1058  &   25.2    &   38      & 313  &  8.24   \\
  160 & 16.5 &     42   &   968  &   23.1    &   41      & 276  &  6.73   \\
  190 & 16.5 &    ---   &  ---   &   ---     &   39      & 339  &  8.69   \\
  220 & 16.5 &    ---   &  ---   &   ---     &   40      & 433  & 10.82   \\
  250 & 16.5 &    ---   &  ---   &   ---     &   47      & 652  & 13.87   \\[1ex]
\hline\hline
   70 & 21.0 &     43   &   585  &  13.60    &   40      & 200  &  5.00  \\
  100 & 21.0 &     42   &  1004  &  23.90    &   42      & 301  &  7.17  \\
  130 & 21.0 &     44   &   914  &  20.77    &   39      & 274  &  7.03  \\
  160 & 21.0 &     25   &   896  &  35.84    &   43      & 507  & 11.79  \\
  190 & 21.0 &    ---   &  ---   &   ---     &   46      & 637  & 13.85  \\
  220 & 21.0 &    ---   &  ---   &   ---     &   41      & 855  & 20.85  \\
  250 & 21.0 &    ---   &  ---   &   ---     &   43      & 891  & 20.72  \\[1ex]
\hline\hline
\end{tabular}
\label{tab:it-arpack-jd}
\end{table}

In the sequel we compare the traditional \cwi {} method with the
\jd\ code \jdbsym \cite{Geu02} using \ilupack\ as preconditioner.
Table \ref{tab:cpu-cwi-jd} shows that switching from \arpack\ to
\jd\ in this case improves the total method by another factor $6$
or greater. For this reason \jd\ together with \ilupack\ will be
used as default solver in the following. The numerical results in
Table \ref{tab:cpu-cwi-jd} show a dramatic improvement in the
computation time by using \ilupack-based \jd. Although this new
method slows down for smaller $w$ due to poorer diagonal
dominance, a gain by orders of magnitude can still be observed.
For $w=16.5$ and larger, even more than three orders of magnitude
in the computation time can be observed. Hence the new method
drastically outperforms the \cwi\ method while the memory
requirement is still moderate.
%
%
%
%
\begin{table}[htb]
    \centering
    \caption{CPU times in seconds and memory requirements in GB to compute five eigenvalues
             closest to $\lambda=0$
             with \cwi {} and \jd\ {}  using \ilupackmatch {} for the shift-and-invert
             technique. '$\ddagger$' indicates that the convergence of the
             method was too slow. For \cwi\ and $M=100$, not all $5$ eigenvalues converged
             successfully, so the eigenvector reconstruction finished quicker,
             leading to variances in the CPU times ($^*$).}
\begin{tabular}{|r|r|c|c|c|c|}\hline
        \hline
   M  &  W   &  \multicolumn{2}{c|}{\cwi} & \multicolumn{2}{c|}{\jd\ {} } \\
      &      &  \multicolumn{2}{c|}{    } & \multicolumn{2}{c|}{\ilupackmatch\ {} } \\
      &      &   time    &   mem. &    time  & mem.\\ \hline\hline
   70 & 12.0 &  20'228   &   0.11 &   1'314  &    0.95 \\
  100 & 12.0 & 148'843   &   0.32 &   8'522  &    2.93 \\
  130 & 12.0 &$\ddagger$ & $\ddagger$        &  56'864  &    8.06 \\
\hline\hline
   70 & 16.5 &  15'100   &   0.11 &     258   &    0.34 \\
  100 & 16.5 &255'842$^*$&   0.32 &     978   &    0.98 \\
  130 & 16.5 &$\ddagger$ & $\ddagger$        &   2'895   &    2.16 \\
  160 & 16.5 &$\ddagger$ & $\ddagger$        &   5'860   &    4.05 \\
  190 & 16.5 &$\ddagger$ & $\ddagger$        &  16'096   &    6.75 \\
  220 & 16.5 &$\ddagger$ & $\ddagger$        &  30'160   &   10.58 \\
  250 & 16.5 &$\ddagger$ & $\ddagger$        &  86'573   &   \\[1ex]
\hline\hline
   70 & 21.0 & 14'371    &   0.11 &     255   &    0.26 \\
  100 & 21.0 &331'514$^*$&   0.32 &   1'134   &    0.75 \\
  130 & 21.0 &$\ddagger$ & $\ddagger$        &   1'565   &    1.65 \\
  160 & 21.0 &$\ddagger$ & $\ddagger$        &   4'878   &    3.08 \\
  190 & 21.0 &$\ddagger$ & $\ddagger$        &  11'032   &    5.22 \\
  220 & 21.0 &$\ddagger$ & $\ddagger$        &  28'334   &         \\
  250 & 21.0 &$\ddagger$ & $\ddagger$        &  48'223   &         \\
\hline\hline
\end{tabular}
\label{tab:cpu-cwi-jd}
\end{table}

One key to the success of the preconditioner is based on the
threshold $\kappa$ which bounds the growth of $L^{-1}$. Already
for a small example  such as $M=70$ significant differences can be
observed. As we show in Table \ref{tab:inverse-bound}, increasing
the bound by a factor $2$ from $\kappa=5$
 up to  $\kappa=10$ and $\kappa=20$ leads to an enormous increase of
fill. Here we measure the fill of the incomplete $LDL^T$
factorization relative to the non-zeros of the original matrix. By
varying the drop tolerance $\varepsilon$ we also see that the
dependence on $\kappa$ is much more significant than the
dependence of $\varepsilon$. Roughly speaking, the ILU
decomposition becomes twice as expensive when $\kappa$ is replaced
by $2\kappa$ and so does the fill-in. The latter is crucial since
memory constraints severely limit the size of the application that
can be computed.
\begin{table}[htb]
    \centering
    \caption{The influence of the inverse bound $\kappa$ on the amount of memory. For $M=70$ compare for different thresholds how the fill-in $nnz(LDL^T)/nnz(A)$ varies depending on $\kappa$ and state the computation time in seconds.}
    \begin{tabular}{|l|rrr|rrr|rrr|}\hline
        \hline
    $\varepsilon$ &  \multicolumn{3}{c|}{$\kappa=5$} &  \multicolumn{3}{c|}{$\kappa=10$} &  \multicolumn{3}{c|}{$\kappa=20$}  \\
      & fill &  time    &  total  &  fill &  time  &  total & fill & time  & total  \\
      &  &  $LDL^T$   &   time  &   &  $LDL^T$   &   time &  &  $LDL^T$ &  time \\\hline\hline
          $0.01$ &   $5.4$ & $3.7\cdot 10^1$ & $8.7\cdot 10^2$ &  $8.7$ & $6.7\cdot 10^1$ & $5.0\cdot 10^2$& $15.2$ & $1.6\cdot 10^2$ & $4.8\cdot 10^2$\\
   $0.005$ &   $6.8$ & $5.4\cdot 10^1$ & $4.4\cdot 10^2$ & $11.0$ & $1.0\cdot 10^2$ & $3.8\cdot 10^2$& $19.1$ & $2.3\cdot 10^2$ & $5.0\cdot 10^2$\\
 $0.0025$ &   $8.6$ & $8.1\cdot 10^1$ & $3.1\cdot 10^2$ & $13.8$ & $1.5\cdot 10^2$ & $3.6\cdot 10^2$& $24.1$ & $3.4\cdot 10^2$ & $6.0\cdot 10^2$\\
          $0.001$ &  $11.7$ & $1.3\cdot 10^2$ & $3.0\cdot 10^2$ & $18.0$ & $2.3\cdot 10^2$ & $4.1\cdot 10^2$& $32.1$ & $5.4\cdot 10^2$ & $7.8\cdot 10^2$\\ \hline\hline
 \end{tabular}
    \label{tab:inverse-bound}
\end{table}

In Table \ref{tab:mem-hb} we show how \jdbsym\ and \ilupackmatch\
perform when instead of periodic boundary conditions, we use hard
wall boundaries, i.e.,
$x_{0;j;k}=x_{i;0;k}=x_{i;j;0}=x_{M+1;j;k}=x_{i;M+1;k}=x_{i;j;M+1}=0$
for all $i,j,k$. This is sometimes of interest in the Anderson
problem and generally, it is expected that for large $M$, the
difference in eigenvalues and eigenvectors becomes small when
compared to the standard periodic boundaries. In addition, we also
show results for the so-called off-diagonal Anderson problem
\cite{CaiRS99}. Here, we shift the diagonal to a constant
$\sigma=1.28$ and incorporate the randomness by setting the
off-diagonal elements of $A$ to be uniformly distributed in
$[-1/2, 1/2]$. The graph of the matrix $A$ remains the same. These
values correspond --- similarly to $w_c= 16.5$ used before for
purely diagonal randomness --- to the physically most interesting
localization transition in this model \cite{CaiRS99}. We note that
using hard wall boundary conditions instead of periodic boundary
conditions leads to slightly less fill but increases the number of
iteration steps as can be seen in Table \ref{tab:mem-hb}. This
conclusions carries over to the off-diagonal Anderson problem,
where the memory consumption is less but the iterative part takes
even longer. In principle our results could be improved if we were
to switch to a smaller threshold $\varepsilon$ than the here
uniformly applied $\varepsilon=1/sqrt(N)$.

\begin{table}[htb]
    \centering
    \caption{Difference in performance for our standard problem with
      periodic boundary conditions, the problem with hard wall
      conditions, and the inverse problem with random numerical entries
      in the off-diagonal elements. Memory requirement (in GB) and CPU times
      (in seconds) to compute at the transition the eigenvectors corresponding to the five
      eigenvalues closest to $\lambda=0$ with shift-and-invert \jd{} and the
      \ilupackmatch{} solver using symmetric weighted matchings.}
    \begin{tabular}{|r|c|c|c|c|c|c|}\hline
        \hline
    N &  \multicolumn{2}{c|}{Periodic} &\multicolumn{2}{c|}{Hard
    wall}&\multicolumn{2}{c|}{Inverse} \\
      & time    &  memory    &  time  & memory & time   & memory \\\hline\hline
   70 &    258  &    0.34    &    282 &  0.31  &    457 & 0.27   \\
  100 &    980  &    0.98    &    969 &  0.94  &  2'075 & 0.79   \\
  130 &  5'244  &    2.16    &  2'090 &  2.07  &  6'472 & 1.76   \\
  160 &  9'958  &    4.05    &  5'661 &  3.90  & 11'975 & 3.30   \\
  190 & 14'742  &    6.75    & 13'431 &  6.62  & 27'488 &        \\ \hline\hline
 \end{tabular}
    \label{tab:mem-hb}
\end{table}
\section{Conclusion}
\label{sect:concl}
We have shown that modern approaches to preconditioning based on
symmetric matchings and multilevel preconditioning methods lead to
an astonishing increase in performance and available system sizes
for the Anderson model of localization. This approach is not only
several orders of magnitudes faster than the traditional \cwi {}
approach, it also consumes only a moderate amount of memory thus
allowing to study the Anderson eigenproblem for significantly
larger scales than ever before.

Let us briefly recall the main ingredients necessary for this
progress: At the heart of the new approach lies the use of
symmetric matchings \cite{HagS04} in the preconditioning stage of
the inverse-based incomplete factorization preconditioning
iterative method \cite{ilupack2.0}. Furthermore, the
preconditioning itself is of a multi-level type, complementary to
the often used full-pivoting strategies. Next, the inverse-based
approach is also of paramount importance to keep the fill-in at a
manageable level (see Table \ref{tab:inverse-bound}). And last, we
emphasize that these results, of course, reflect our selected
problem class: to compute a few of the interior eigenvalues and
associated eigenvectors for a highly indefinite symmetric matrix
defined by the Anderson model of localization.

The performance increase by several orders of magnitude (see Table
\ref{tab:cpu-cwi-jd}) is solely due to our use of new and improved
algorithms. Combined with advances in the performance to cost
ratio of computing hardware during the last 6 years period,
current preconditions methods makes it possible to solve those
problems quickly and easily which have been considered by far too
large until recently \cite{CaiMRS02}. Even for $N\times N$
matrices as large as $N = 64 \cdot 10^6$, it is now possible
within a few days to compute the interior eigenstates of the
Anderson problem.

The success of this method indicates that it might also be
successfully applied to other large-scale problems arising in
(quantum) physics.

\section*{Acknowledgement}
We gratefully acknowledge discussions and implementation help from
A.\ Croy and C.\ Sohrmann in the initial stages of the project.



\end{document}